\documentclass[a4paper,twoside,fleqn,british]{elsarticle}
\usepackage{prettyref}
\usepackage{float}
\usepackage{units}
\usepackage{mathrsfs}
\usepackage{amsmath}
\usepackage{amssymb}
\usepackage{graphicx}
\usepackage{esint}
\usepackage{nameref}
\biboptions{sort&compress}
\usepackage{babel}
\begin{document}
\begin{frontmatter}{}\par 
	\title{A statistical interpretation of biologically inspired growth models}\par \tnotetext[t1]{This document is a collaborative effort.}\par 
\author{A. Samoletov\corref{cor1}\fnmark[1]}\par \ead{A.Samoletov@liverpool.ac.uk}\par 
\author{B. Vasiev\corref{cor2}}\par \ead{B.Vasiev@liverpool.ac.uk}\par \cortext[cor1]{Corresponding author}\par 
\cortext[cor2]{Corresponding author}\par 
\fntext[1]{Also at Galkin Institute for Physics and Engineering}\par 
\address{The University of Liverpool, Liverpool, UK}
\begin{abstract} Biological entities are inherently dynamic. As such,
various ecological disciplines use mathematical models to describe
temporal evolution. Typically, growth curves are modelled as sigmoids,
with the evolution modelled by ordinary differential equations. Among
the various sigmoid models, the logistic and Gompertz equations are
well-established and widely used in fitting growth data in the fields
of biology and ecology. This paper suggests a statistical interpretation
of the logistic equation within the general framework. This interpretation
is based on modelling the population--environment relationship, the
mathematical theory of which we discuss in detail. By applying this
theory, we obtain stochastic evolutionary equations, for which the
logistic equation is a limiting case. The prospect of modifying logistic
population growth is discussed. We support our statistical interpretation
of population growth dynamics with test numerical simulations. We
show that the Gompertz equation and other related models can be treated
in a similar way. 
\end{abstract} 
\par \end{frontmatter}{}

\section{Introduction }

Dynamics of a single species population is a prototype for mathematical
modelling in ecology and is commonly described by the following differential
equation, 
\begin{equation}
\dot{x}=xf(x),\label{eq:Kolm}
\end{equation}
where $x$ is the population density and $f(x)$ describes the\emph{
per capita} growth. The observed upper limit for population growth
determines the level of saturation known as carrying capacity. The
simplest form of the associated dynamics, the Verhulst-Pearl \citep{Verhulst1838,pearl1924curve,kingsland1982}
logistic equation, 
\begin{equation}
\dot{x}=rx\left(1-\frac{x}{K}\right),\label{eq:logistic}
\end{equation}
where $r$ is a time scale so that $rt$ is dimensionless time and
$K$ is the carrying capacity, is often and successfully used to model
population growth \citep{kot2001elements,Murray2002}. The modelling
of population dynamics using equation \eqref{eq:Kolm} is based on
a special choice of mathematical expressions for the laws governing
population growth, that is, the function $f(x)$, so that the environment
is treated as a static reservoir. The logistic equation \eqref{eq:logistic}
inherits the requested generic properties of the equation \eqref{eq:Kolm}.
To substantiate a choice $f(x)=r\left(1-K^{-1}x\right)$, the data
of experimental observations and heuristic reasoning are usually used
\citep{Verhulst1838,pearl1924curve,kingsland1982,kot2001elements,Murray2002}.

Another approach to justifying the equation \eqref{eq:logistic} was
initiated by Volterra \citep{volterra1939calculus} and was based
on finding a certain minimum principle leading to the logistic law
of population growth, that is, to the equivalent Euler-Lagrange equation.
The problem is to find the Lagrangian that has the required form.
This approach was discussed in the literature and several forms of
the Lagrangian were proposed \citep{volterra1939calculus,leitmann1972minimum,gatto1988functional,webb1995hamilton,pawlowski2006dynamic}.
Also some speculations were made regarding the universality of the
functional approach in the context of ecological problems \citep{webb1995hamilton,pawlowski2006dynamic,wilhelm2000goal}.

Philosophically, one can trivially assert that all subsystems of the
ecosphere are in interaction and interdependence. However, this statement
is meaningless, since only a relatively small number of characteristic
variables need to be taken into account in specific models. Since
reliable information about the actual state of the entire ecosphere
is unknown, its impact on small ecosystems can only be accounted for
phenomenologically, based on plausible reasoning.

A statistical approach to population dynamics, which can be either
deterministic or stochastic, is as follows. Consider an ensemble of
identical populations, differing only in initial size, together with
the corresponding phase space. Assuming that the ensemble is initially
characterised by a probability density on the phase space, we interpret
the population dynamics as an evolution of this ensemble, and we are
interested in the asymptotic density at large times. Traditional deterministic
population dynamics can be represented by a stable asymptotic behaviour
\citep{may2019stability}, so that the resulting densities are singular.
If the dynamics can be perturbed in some way, we must consider smooth
density functions. In the context of dynamical systems theory we are
interested in invariant densities. Ergodicity means equivalence of
ensemble averaging and time averaging. This property is the basis
of statistical approach to dynamics, and therefore the ergodic hypothesis
is necessary for our analysis.

By moving towards a mathematical formulation of a statistical approach
to population dynamics, we recognise that natural populations do not
exist and grow in isolation, but are in dynamic contact with the environment.
It is clear that population evolution is influenced by changes in
static (or quasi--static) environmental conditions, but to better
understand population growth, dynamic models must also account for
population--environment interactions, environmental responses to
a growing population, and internal population and environmental processes
that may be only partially known, introducing a degree of uncertainty
into population dynamics. Indeed, any population placed in the environment
must perturb it to some extent and will itself be subject to the backward
influence of this perturbation. Depending on the time scale of the
population growth the environment can be divided into two parts: the
part involved in the joint dynamics with the population and the unperturbed
part, which determines the general statistical properties of the population--environment
coevolution and the long-term homeostasis of the population, if it
exists.

To ensure an equilibrium state, population growth must be controlled
by resource balance, population abundance and a variety of other mechanisms.
The dynamic processes of fluctuation and relaxation must be defined
accordingly. However, we do not specify explicit mechanisms. Instead,
we define general statistical hypotheses that cover many different
biological processes. Specific biological processes need to be considered
separately.

In the context of population dynamics and ecosystem evolution, the
statistical approach, as we formulate it, contains the possibility
of implementing different scenarios of transition from the initial
state to the asymptotic state of statistical equilibrium, which implies
qualitatively different dynamic scenarios. In our case, the goal is
the statistical interpretation of a very special process - the logistic
law of population growth.

To present a statistical interpretation of the logistic equation,
that is, a rather simple but important mathematical model of population
growth, we begin with a brief statement of the relevant mathematical
assumptions that combine the dynamic principle for ensemble control
tools, the definition of the environmental temperature and the corresponding
invariant density, the ergodic hypothesis and methods of stochastic
analysis \citep{SamoletovVasiev2017,samoletov2021advanced,SamoletovDettmannChaplain2007,SamoletovDettmannChaplain2010}.
Then, based on plausible probabilistic reasoning, we formulate and
explore the statistical interpretation of the logistic equation.

\section{Mathematical formalism\label{sec:Mathematical-formalism}}

In this section, we briefly summarise the mathematical concepts necessary
to understand the statistical interpretation of the logistic equation.

Let a population $\mathscr{P}$ be placed in an environment $\mathscr{E}$,
a system of large (infinite) number of phase variables, that determines
the general statistical properties of the entire ecosystem. The population
has to disturb the environment to some extent, and will be affected
by this disturbance. Let the environment $\mathscr{E}$ be divided
into the part participating in the joint dynamics with the population,
$\mathscr{E}^{*}$, and the rest of the ecological system, $\mathscr{E}\setminus\mathscr{E^{*}}$,
which determines statistical properties of population and environment
coevolution, and long-term population homeostasis. Subsystems $\mathscr{P}$
and $\mathscr{E^{*}}$ interacting with the global environmental reservoir
$\mathscr{E}\setminus\mathscr{E^{*}}$ can be fluctuating, while $\mathscr{E}\setminus\mathscr{E^{*}}$
remains unchanged, determining the general statistical properties
of the whole system, $\mathscr{P}+\mathscr{E^{*}}$. Let us make an
important assumption that the population $\mathscr{P}$ and its environment
$\mathscr{E}^{*}$ participating in the joint dynamics are statistically
independent in equilibrium. Of course, the actual description of $\mathscr{E}^{*}$
cannot be done in advance and depends on the population, the environment
and the experiments used to extract information on population growth,
as they determine the temporal and spatial scales of measurements
and the corresponding interpretation of data. Thus, the system $\mathscr{E}^{*}$
is not predetermined, but depends on a number of factors that do not
affect the asymptotic statistical behaviour of the population. Among
these factors, the possibility of various evolutionary transients
should be emphasised. The study of transients in population dynamics
and ecology is an intensively studied problem \citep{hastings2004transients,Hastings2018}.
Practically, to describe $\mathscr{E}^{*}$, one can try one or another
set of variables to construct consistent coevolution equations and
investigate their properties.

To consider the interaction and joint evolution of systems $\mathscr{P}$
and $\mathscr{E}^{*}$ as described above, we first need to define
the dynamical system, $\mathrm{\mathscr{S}}=\left(\mathcal{M}^{+},G(z)\right)$,
which is a direct product of noninteracting (autonomous) systems $\mathscr{P}=\left(\mathcal{M},g(x)\right)$
and $\mathscr{E}^{*}=\left(\mathcal{M}^{*},g^{*}(y)\right)$, that
is, $\mathscr{S}=\mathscr{P}\times\mathscr{E}^{*}=\left(\mathcal{M}\oplus\mathcal{M}^{*},g(x)\times g^{*}(y)\right)$,
where $\mathcal{M}$ and $\mathcal{M}^{*}$ are phase spaces of systems
$\mathscr{P}$ and $\mathscr{E}^{*}$ correspondingly. In other words,
we consider a simple combination of two independent systems into one
so that $z=\left(x,y\right)\in\mathcal{M}^{+}=\mathcal{M}\oplus\mathcal{M}^{*}$
and $\dot{z}=G(z)$, where $G(z)=g(x)\times g^{*}(y)$. When $\mathscr{P}$
and $\mathscr{E}^{*}$ are considered as systems involved in joint
evolution, such a separation into noninteracting systems becomes impossible
and we have to consider dynamics in the general form, 
\begin{equation}
\dot{z}=G^{+}(z).\label{eq:Eq-z}
\end{equation}
However, it is important that, as in the case of noninteracting systems,
the invariant density $\sigma^{+}$ for the combined system $\mathscr{S}$
is the direct product of the invariant densities $\sigma$ and $\sigma^{*}$
for the systems $\mathscr{P}$ and $\mathscr{E}^{*}$, 
\begin{equation}
\sigma^{+}\left(z\right)=\sigma\left(x\right)\times\sigma^{*}\left(y\right),\label{eq:Denc-z}
\end{equation}
that is, the systems $\mathscr{P}$ and $\mathscr{E}^{*}$ are statistically
independent in the equilibrium state, provided that the equilibrium
state of the entire ecosystem exists, where $\mathscr{P}$ and $\mathscr{E}^{*}$
are small subsystems.

As preparation for what follows, we briefly summarise two important
concepts: (1) ecological temperature and (2) the dynamic principle
for statistical ensemble control tools.

\subsection*{Ecological temperature}

The concept of ecological temperature is based on the concept of temperature
expression described in \citep{SamoletovVasiev2017,samoletov2021advanced}.
Let us briefly outline the relevant details.

Let the probability density $\sigma\left(x\right),\:x\in\mathcal{M}$
be given. Define the function $h(x):\mathcal{M}\rightarrow\mathbb{R}$,
$h(x)\propto-\vartheta\ln\sigma\left(x\right)$, where $\vartheta>0$
is a parameter, so that $h(x)$ is a sufficiently smooth function,
bounded from below and growing at infinity, $h(x)\geq a\left|x\right|^{b}$
for some $a>0,\:b>0$, that is, a coercive function. We can now express
the equilibrium probability density function in the form
\begin{gather}
\sigma_{\vartheta}\left(x\right)\propto\exp\left\{ -\vartheta^{-1}h(x)\right\} ,\quad x\in\mathcal{M}.\label{eq:Denc-x}
\end{gather}
The probability density $\sigma_{\vartheta}^{*}\left(y\right)$ is
understood in the same way. That is, assume that the equilibrium probability
density function associated with the environment is $\sigma^{*}\left(y\right),\,y\in\mathcal{M}^{*}$.
Then we define the function $h^{*}(y):\mathcal{M}^{*}\rightarrow\mathbb{R}$,
$h^{*}(y)\propto-\vartheta\ln\sigma^{*}\left(y\right)$ and represent
the density function in the form 
\begin{gather}
\sigma_{\vartheta}^{*}\left(y\right)\propto\exp\left\{ -\vartheta^{-1}h^{*}(y)\right\} ,\quad y\in\mathcal{M}^{*}.
\end{gather}

The function $\Theta(x,\vartheta)$, $\Theta:\mathcal{M}\times\mathbb{R}_{+}\rightarrow\mathbb{R}$
is called an ecological temperature expression (abbreviated as $\vartheta$-expression)
if it explicitly depends on the parameter $\vartheta$ and satisfies
the conditions, 
\begin{equation}
\mathbb{E}_{\vartheta}\left\{ \Theta(x,\vartheta)\right\} =0\quad\textrm{for all}\quad\vartheta>0,\label{eq:TEx-1}
\end{equation}
where $\mathbb{E}_{\vartheta}\left\{ \ldots\right\} =\int_{\mathcal{M}}(\ldots)d\mu_{\vartheta}(x)$
is the mathematical expectation, and $d\mu_{\vartheta}(x)=\sigma_{\vartheta}(x)dx$
is the probability distribution. The $\vartheta$-expression \eqref{eq:TEx-1}
is defined up to a constant factor, possibly depending on $\vartheta$,
and an additive function $\psi(x)$ such that $\int_{\mathcal{M}}\psi(x)d\mu_{\vartheta}(x)=0$.
We consider $\Theta(x,\vartheta)$ as an analytic function of the
real parameter $\vartheta$, $\Theta(x,\vartheta)=\Theta_{0}(x)+\Theta_{1}(x)\vartheta+\ldots$,
whose first term $\Theta_{0}(x)$ has a nonzero expectation $\mathbb{E}_{\vartheta}\left\{ \Theta_{0}(x)\right\} \neq0$.
Usually $\Theta(x,\vartheta)$ is a polynomial in $\vartheta$ \citep{samoletov2021advanced}.
This is the context in which definition \eqref{eq:TEx-1} should be
understood. In practice, although $\psi(x)$ functions are not true
$\vartheta$-expressions, they can nevertheless be used to generate
deterministic equations of motion consistent with the dynamic principle,
as defined below. However, in this case the required property of ergodicity
seems to be more than doubtful (e.g. \citep{legoll2009non}).

Simply put, the specific meaning of the definition of $\vartheta$--expression
is that it allows us to express the value of the ecological temperature
$\vartheta$ as an average of some function of the dynamic variables
of either the population or the environment.

We consider $\mathscr{P}$ and $\mathscr{E}^{*}$ as parts of a large,
infinitely large, ecosystem that determines the equilibrium statistical
properties of both the population $\mathscr{P}$ and environment $\mathscr{E}^{*}$,
provided that such an equilibrium state exists. Thus, it should be
assumed that the ecological temperatures of the population $\mathscr{P}$
and the environment $\mathscr{E}^{*}$ coincide (see \citep{khinchin1949mathematical}
to argue for this). This equality of the parameter $\vartheta$ has
direct implications for the mathematical formulation of the theory.

The set of all $\vartheta$-expressions for an arbitrary but fixed
value of the parameter $\vartheta$ is a linear system in which the
operations of addition and multiplication by numbers are defined in
the usual way. In order to study and use the properties of $\vartheta$-expressions,
it will be necessary to interpret them as elements of either the space
$L_{1}$ (summable expressions) or $L_{2}$ (square summable expressions).
Such an interpretation is appropriate here, since the $\vartheta$-expressions
we are considering are bounded from below and grow at infinity no
faster than a polynomial.

For the combined system $\mathscr{S}$ a temperature expression satisfies
the condition, 
\[
\int_{\mathcal{M}^{+}}\Theta^{+}(z,\vartheta)d\mu_{\vartheta}^{+}(z)=\int_{\mathcal{M}^{+}}\Theta^{+}(x,y,\vartheta)d\mu_{\vartheta}(x)d\mu_{\vartheta}^{*}(y)=0\quad\textrm{for all}\quad\vartheta>0,
\]
and, provided that $\int_{\mathcal{M}}\left|\Theta^{+}(z,\vartheta)\right|d\mu_{\vartheta}(z)<\infty$,
it follows from Fubini's theorem that 
\[
\int_{\mathcal{M}}\Theta^{+}(z,\vartheta)d\mu_{\vartheta}(x)=\Theta^{*}(y,\vartheta)\quad\mathrm{and}\quad\int_{\mathcal{M}^{*}}\Theta^{+}(z,\vartheta)d\mu_{\vartheta}^{*}(y)=\Theta(x,\vartheta)
\]
are $\vartheta$-expressions as defined above. For more information
about the properties and selection of $\vartheta$--expressions,
see \citep{samoletov2021advanced}.

\subsection*{Dynamic principle}

The dynamic principle \citep{SamoletovVasiev2017,samoletov2021advanced}
for development of statistical ensemble control tools is based on
the assumption of ergodicity, that is, the averaging can equally be
interpreted either as an ensemble average or as a time average for
a single trajectory. To unify the notation, we denote the result of
the averaging by $"\sim"$ and write, for example, in relation to
the $\vartheta$-expression, $\Theta(x,\vartheta)\sim0$, assuming
that $\underset{t\rightarrow\infty}{\lim}\frac{1}{t}\intop_{0}^{t}\Theta(x(t'),\vartheta)dt'=0$.

In the statistical description of a dynamical system $\dot{x}=g(x)$,
the concept of the first integral plays a central role \citep{khinchin1949mathematical}.
The function $h\left(x\right)$ is the first integral if and only
if $\nabla h\left(x\right)\cdot g\left(x\right)=0$ for all $x\in\mathcal{M}$.
When populations are involved in coevolution with the environment,
then $h\left(x\right)$ is no longer the first integral. This is also
true when the effects of the ecosystem on the population are treated
as random perturbations. To describe such coevolution with a given
invariant measure, the dynamic principle for ensemble control tools
is used, which leads to consistent dynamic equations. The invariant
measure can be \emph{a priori} or can be derived from experimental
data. 

Let us assume that the equations of motion are of the form \eqref{eq:Eq-z}
(or $\dot{z}=G_{\omega}^{+}(z)$ in the case of stochastic dynamics,
indicated by the sub-index $\omega$) and the invariant probability
density is of the form \eqref{eq:Denc-z}, explicitly $\sigma_{\vartheta}^{+}\left(z\right)\propto\exp\left\{ -\vartheta^{-1}h^{+}\left(z\right)\right\} $,
where $h^{+}\left(z\right)=h\left(x\right)+h^{*}\left(y\right)$.
If we denote either $\Gamma\left(z\right)\equiv\nabla_{z}h^{+}(z)\cdot G^{+}(z)$
for deterministic, either $\Gamma\left(z\right)\equiv\mathbb{E}_{\omega}\left\{ \nabla_{z}h^{+}(z)\cdot G_{\omega}^{+}(z)\right\} $
for stochastic dynamics, where $\mathbb{E}_{\omega}\left\{ \ldots\right\} $
denotes the averaging over all realisations of random processes, then
the dynamic principle postulates the following functional relation:
\begin{equation}
\Gamma\left(z\right)\propto\Theta^{+}(z,\vartheta)\sim0.\label{eq:DP-1}
\end{equation}
Solutions to the functional equation \eqref{eq:DP-1} represent possible
transient population evolution scenarios compatible with a given invariant
density. To find a particular solution to the equation \eqref{eq:DP-1},
it is necessary to specify an invariant density and to select an admissible
$\vartheta$--expression corresponding to that density. The criterion
for such a selection is determined by the nature of the problem to
be solved, for example, it may be the simplest admissible expression,
\emph{i.e.}, a first-order polynomial in $\vartheta$. However, the
range of possible selections is wide \citep{samoletov2021advanced,SamoletovVasiev2017,SamoletovDettmannChaplain2007,SamoletovDettmannChaplain2010}.

To understand the practical value of this rather abstract mathematical
scheme, let us start simple and consider a conceptual example. To
this end, we will consider a statistical interpretation of the logistic
equation.

\section{Statistical interpretation of the logistic equation}

\subsection{Preliminaries}

Consider $\mathscr{P}$ as a single species homogeneous population
with density $x\in\mathbb{R}_{+}$, placed in an environment such
that a state of equilibrium exists. This is not a static equilibrium.
Population growth depends on available resources and other environmental
conditions, as well as the size of the population itself, to ensure
a dynamic statistical equilibrium. These processes, which involve
the environment in coevolution, must have balanced dynamics where
the fluctuation and relaxation processes must be appropriately specified.
To formulate a statistical approach to population dynamics, we will
assume that the state of equilibrium and the corresponding ecological
temperature $\vartheta$ are determined by the ecosystem as a whole,
of which the population is a small part.

According to the theoretical scheme presented in Section \ref{sec:Mathematical-formalism},
we must first establish the invariant density for the system $\mathscr{S}=\mathscr{P}+\mathscr{E}^{*}$,
\emph{i.e.} $\sigma_{\vartheta}^{+}\left(z\right)=\sigma_{\vartheta}\left(x\right)\times\sigma_{\vartheta}^{*}\left(y\right)$.

\subsubsection{\emph{The environment\label{subsec:env}}}

Although it is not compulsory, let us suppose that the environment
$\mathscr{E^{*}}$ is characterised by a variable $y\in\mathbb{R}$
with a probability density $\sigma_{\vartheta}^{*}(y)$. By necessity,
$y$ is a collective variable that incorporates a number of environmental
factors and processes. Thus, it is reasonable to guess the Gaussian
statistics, that is, 
\begin{equation}
\sigma_{\vartheta}^{*}(y)\propto\exp\left\{ -\vartheta^{-1}\frac{1}{2}y^{2}\right\} ,\label{eq:Gauss-y}
\end{equation}
where $y$ is a dimensionless variable. The ecological temperature
$\vartheta$ defines the intensity of the environmental fluctuations.
The Gaussian statistics is completely characterised by the first two
cumulants: $\mathbb{E}_{\vartheta}\left(y\right)=0$, $\mathrm{\mathbb{E}}_{\vartheta}\left\{ \left[y-\mathrm{\mathbb{E}}_{\vartheta}\left\{ y\right\} \right]^{2}\right\} =\mathrm{\mathbb{E}}_{\vartheta}\left\{ y^{2}\right\} =\vartheta$.
Therefore, a linear combination of $y$ and $\left(y^{2}-\vartheta\right)$,
say $c_{1}y+c_{2}\left(y^{2}-\vartheta\right)$, is a $\vartheta$--expression.
It should be noted that (1) parameters $c_{1}$ and $c_{2}$ are allowed
to be functions of the population density variable $x$, (2) even
in the case of Gaussian statistics, 
\begin{equation}
\Theta^{*}(y,\vartheta)=c_{1}y+c_{2}\left(y^{2}-\vartheta\right)\label{eq:ecotemp}
\end{equation}
 will be the simplest $\vartheta$--expression (first-order polynomial
in $\vartheta$) associated with the environment. Other higher-order
polynomial expressions in $\vartheta$ exist {[}12{]}. As an example,
the Chebyshev-Hermite polynomials with the parameter $\vartheta$,
$He_{n}\left(y;\vartheta\right),\:n\in\mathbb{N}_{0}$,
\[
He_{n}(y;\vartheta)=\exp\left\{ \frac{y^{2}}{2\vartheta}\right\} \left(-\vartheta\right)^{n}\frac{d^{n}}{dy^{n}}\exp\left\{ -\frac{y^{2}}{2\vartheta}\right\} ,\quad n\in\mathbb{N}_{0},
\]
(the formula, commonly referred to as the Rodrigues formula) are $\vartheta$--expressions
for the Gaussian probability density \eqref{eq:Gauss-y} for $n\geq1$,
that is, $\mathrm{\mathbb{E}}_{\vartheta}\left\{ He_{n\geq1}(y;\vartheta)\right\} =0$.
Explicitly, $He_{1}(y;\vartheta)=y$, $He_{2}(y;\vartheta)=y^{2}-\vartheta$,
and so forth.

We will explore $\vartheta$--expression \eqref{eq:ecotemp}, which
comprises constant coefficients $c_{1}$ and $c_{2}$, to define ecological
temperature in the simplest way possible (see {Appendix} for details).

\subsubsection{\emph{The population\label{subsec:pop}}}

Let us accept the invariant probability density $\sigma_{\vartheta}\left(x\right)$
in the form \eqref{eq:Denc-x}, where the coercive function $h(x)$
is a subject to define. The $\vartheta$ parameter is the same as
in $\sigma_{\vartheta}^{*}(y)$. The simplest $\vartheta$-expression
associated with the population, that is, a first-order polynomial
in $\vartheta$ with $\vartheta$ as the additive term, has the form,
\begin{equation}
\Theta(x,\vartheta)=xh'(x)-\vartheta,\label{eq:THETA-1}
\end{equation}
where prime denotes derivative (Lagrange notation). The proof is by
direct calculation, that is, 
\[
\mathrm{\mathbb{E}}_{\vartheta}\left\{ xh'(x)\right\} =\left(\intop\exp\left\{ -\vartheta^{-1}h(x)\right\} dx\right)^{-1}\intop xh'(x)\exp\left\{ -\vartheta^{-1}h(x)\right\} dx=\vartheta.
\]
From the Rodrigues type formula, 
\[
\Theta_{n}(x,\vartheta)=\exp\left\{ \vartheta^{-1}h(x)\right\} \left(-\vartheta\right)^{n}\frac{d^{n}}{dx^{n}}\left[\varphi(x,\vartheta)\exp\left\{ -\vartheta^{-1}h(x)\right\} \right],
\]
we can obtain a series of $\vartheta$--expressions, depending on
the choice of the function $\varphi(x,\vartheta)$ (the dependence
on $\vartheta$ is optional), which does not grow faster than a polynomial.
For example,
\[
\Theta_{1}(x,\vartheta)=\varphi(x)h'(x)-\varphi'(x)\vartheta.
\]
If we set $\varphi(x)=x$ then for we get the expression \eqref{eq:THETA-1}.

We set the ecological temperature value by the $\theta$--expression
associated with the environment \eqref{eq:ecotemp}, which seems quite
reasonable. A $\theta$--expression associated with the population
is then generated during the solution of the dynamic principle equation
(further details can be found in the Appendix.).

\subsection{Logistic equation}

To proceed to the statistical interpretation of the logistic equation,
it is necessary to solve the functional equation \eqref{eq:DP-1}
that includes at least the $\vartheta$--expression associated with
the environment, that is, the equation 
\begin{equation}
\mathbb{E}_{\omega}\left\{ h'\left(x\right)g(x,y)+yg^{*}(x,y)\right\} =\Theta^{*}(y,\vartheta)\sim0,\label{eq:DPlog}
\end{equation}
where the $\vartheta$--expression $\Theta^{*}(y,\vartheta)$ is
given by equation \eqref{eq:ecotemp} with nonzero constant coefficients
and the equations of motion are presented in the following form:
\begin{eqnarray*}
\dot{x} & = & g(x,y),\\
\dot{y} & = & g^{*}(x,y).
\end{eqnarray*}
If $c_{2}\neq0$, then among the solutions of equation \eqref{eq:DPlog}
there are only stochastic equations \citep{SamoletovVasiev2017}.
The procedure of determining the functions $g(x,y)$ and $g^{*}(x,y)$
shares similarities with identifying the first integrals for a given
dynamical system. The task of discovering these integrals can be challenging.
However, it is relatively easy to confirm whether a function satisfies
the requirement of a first integral. To derive a particular solution
to the functional equation \eqref{eq:DPlog}, we use the method previously
discussed in another context \citep{SamoletovVasiev2017}. To avoid
excessive complexity of the current argument, the corresponding calculations
are provided in \nameref{sec:Appendix-A}. This section confirms that
the derived equations of motion meet all the required conditions. 

Let us consider the equations of motion:
\begin{flalign}
\dot{x} & =\lambda xy,\nonumber \\
\dot{y} & =-\lambda\Theta(x,\vartheta)-\gamma y+\sqrt{2\gamma\vartheta}\xi(t),\label{eq:StochPopDyn}
\end{flalign}
where $\xi(t)$ is the standard Gaussian white noise, $\mathbb{E}_{\omega}\left\{ \xi(t)\right\} =0$,
$\:\mathbb{E}_{\omega}\left\{ \xi(t)\xi(t')\right\} =\delta(t-t')$,
$\lambda>0$ and $\gamma>0$ are parameters that actually define two
time scales.

First we check the $\vartheta$--expression for the environment,
$\Theta^{*}(y,\vartheta)$, by direct calculation. That is, by substituting
the expressions \eqref{eq:StochPopDyn} into the equation \eqref{eq:DPlog},
we calculate that $\Theta^{*}(y,\vartheta)=\lambda\vartheta y-\gamma\left(y^{2}-\vartheta\right)$
(see Appendix for details of the calculations).

Then, we prove that the density, 
\begin{equation}
	\sigma_{\vartheta}^{+}(x,y)\propto\exp\left\{ -\vartheta^{-1}h(x)\right\} \times\exp\left\{ -\vartheta^{-1}\frac{1}{2}y^{2}\right\} ,\label{eq:probdensity}
\end{equation}
is invariant for dynamics \eqref{eq:StochPopDyn}. Indeed, the Fokker-Planck
equation corresponding to stochastic differential equation \eqref{eq:StochPopDyn}
(\emph{e.g.}, \citep{gardiner2009stochastic,klyatskin2005dynamics})
has the form $\partial_{t}\sigma=\boldsymbol{\mathcal{F}}^{*}\sigma$,
where 
\begin{gather}
	\boldsymbol{\mathcal{F}}^{*}\sigma=-\frac{\partial}{\partial x}\left[\lambda xy\sigma\right]-\frac{\partial}{\partial y}\left\{ \left[-\lambda\left[xh'(x)-\vartheta\right]-\gamma y-\gamma\vartheta\frac{\partial}{\partial y}\right]\sigma\right\} ,\label{eq:FPE}
\end{gather}
is the Fokker-Planck operator. The most straightforward way of writing
down the Fokker-Planck equation \eqref{eq:FPE} is to use the approach
\citep{klyatskin2005dynamics}.

We prove the identity $\boldsymbol{\mathcal{F}}^{*}\sigma_{\vartheta}^{+}(x,y)=0$
by direct calculation, implying that $\sigma_{\vartheta}^{+}(x,y)$
\eqref{eq:probdensity} is the invariant density for stochastic dynamics
\eqref{eq:StochPopDyn}. Explicitly,
\begin{eqnarray*}
-\frac{\partial}{\partial x}\left[\lambda xy\sigma_{\vartheta}^{+}(x,y)\right]-\frac{\partial}{\partial y}\left\{ \left[-\lambda\left[xh'(x)-\vartheta\right]-\gamma y-\gamma\vartheta\frac{\partial}{\partial y}\right]\sigma_{\vartheta}^{+}(x,y)\right\} \\
=-\frac{\partial}{\partial x}\left[\lambda xy\sigma_{\vartheta}^{+}(x,y)\right]-\frac{\partial}{\partial y}\left\{ -\lambda\left[xh'(x)-\vartheta\right]\sigma_{\vartheta}^{+}(x,y)\right\} \\
+\gamma\frac{\partial}{\partial y}\left\{ y\sigma_{\vartheta}^{+}(x,y)+\vartheta\frac{\partial}{\partial y}\sigma_{\vartheta}^{+}(x,y)\right\} \\
=\left\{ -\lambda y+\lambda xyh'(x)\frac{1}{\vartheta}-\lambda\left[xh'(x)-\vartheta\right]y\frac{1}{\vartheta}\right\} \sigma_{\vartheta}^{+}(x,y)\\
+\gamma\frac{\partial}{\partial y}\left\{ \left[y-y\right]\sigma_{\vartheta}^{+}(x,y)\right\}  & \equiv0.
\end{eqnarray*}
One would expect that for the stochastic evolution equation \eqref{eq:StochPopDyn}
the dynamics would be ergodic.

We have arrived at fairly simple population--environment coevolution
equations \prettyref{eq:StochPopDyn}. We can now ask how these dynamic
equations relate to the conventional population growth equations \eqref{eq:Kolm},
in particular the logistic equation \prettyref{eq:logistic}. To answer
this question, consider the limiting case, $\gamma\gg1$ (relaxation
processes are extremely fast) and $\vartheta\rightarrow+0$ (the environment
is in static equilibrium). In this procedure we follow an analogy
with the Kramers problem \citep{kramers1940,Samoletov1999}. As a
result, we arrive at deterministic dynamics involving only the population
variable. Note that the corresponding mathematically consistent passage
to the limit is not trivial, but intuitively the result seems quite
clear. It is worth noting, however, that following the mathematical
scheme of paper \citep{Samoletov1999} we obtain not only the expressions
given here, but also corrections to them of the next order of magnitude,
which allow us to formulate mathematically the conditions of the limit
transition, including the conditions on the smoothness of the function
$h(x)$ \citep{Samoletov1999}. In the context of our problem, there
appears to be no need for special attention to be paid to such mathematical
details. For more mathematical details, please refer to \citep{mel1991kramers}
and \citep{Samoletov1999}, and the literature cited in them.

Thus, passing to the limit $\gamma\gg1$, $\dot{y}=0$, and $\vartheta\rightarrow+0$,
we get $y=-\gamma^{-1}\lambda\,xh'(x)$, and hence the correspondence
\[
xf(x)=-\lambda^{2}\gamma^{-1}x^{2}h'(x).
\]
For the logistic population growth rate, that is, $f(x)=r\left(1-K^{-1}x\right)$,
setting $\lambda^{2}\gamma^{-1}=r$, we get the following expression
for the function $h(x)$, 
\begin{gather}
h(x)=K^{-1}x-\ln x.\label{eq:log-h}
\end{gather}
We shall say that $h(x)$ \prettyref{eq:log-h} is the logistic $h$--function.
Let us now substitute this $h$--function into the equations \eqref{eq:StochPopDyn}
and thus obtain a stochastic analogue of the logistic equation \eqref{eq:logistic},
which to some extent describes the coevolution of the population and
the environment, 
\begin{alignat}{1}
\dot{x} & =\lambda xy,\nonumber \\
\dot{y} & =-\lambda\left[K^{-1}\left(x-K\right)-\vartheta\right]-\gamma y+\sqrt{2\gamma\vartheta}\xi(t),\label{eq:StatLog}
\end{alignat}
where $\lambda^{2}\gamma^{-1}=r$.

To illustrate the difference between the population evolution described
by the \eqref{eq:logistic} and \eqref{eq:StatLog} equations, we
performed a test numerical simulation of these dynamic equations.
This simulation also allowed us to test the validity of our basic
assumptions. We simulated these equations using the Euler scheme with
a time step $dt=0.001$. In all simulations, we keep $K=1$, $\lambda=1$,
and $\gamma=50$ (the latter to stay close to the logistic equation
\eqref{eq:logistic}), but vary the values of $\vartheta$.

\begin{figure}[h]
\centering{}\includegraphics[scale=0.33]{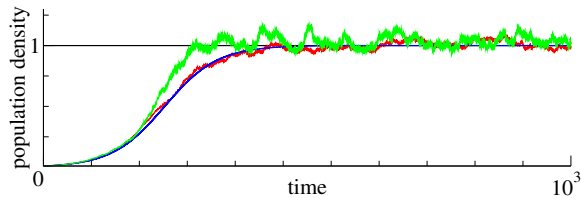}\caption{(colour online) Population density versus time. The various curves
correspond to the Verhulst-Pearl logistic equation \eqref{eq:logistic}
(blue curve) and a stochastic analogue of the logistic equation \eqref{eq:StatLog}
at two values of the $\vartheta$ parameters: $\vartheta=0.001$ (red
curve) and $\vartheta=0.005$ (green curve). Initial values are $x=0.01$,
$y=0$.}
\label{stlog1} 
\end{figure}

\begin{figure}[H]
\centering{}\includegraphics[clip,scale=0.33]{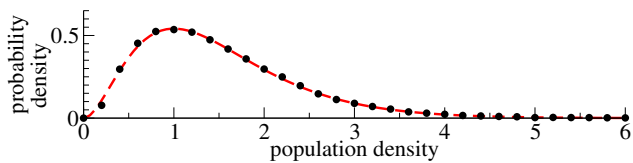}\caption{(colour online) Population probability density versus population density.
The (red) dash lines show the theoretical result given by the logistic
$h$--function. The (black) circles show the normalised histogram
generated at $\vartheta=0.5$ by running numerical simulation for
quite a long time, up to $t_{max}=10^{5}$. Observe, the circles follow
the theoretical curve exactly.}
\label{stlog2} 
\end{figure}

The results of the simulations confirm our assumptions and expectations.
Figure \ref{stlog1} shows that the concept of a statistical approach
to interpreting the logistic equation is reasonable. In the case of
$\gamma=50\gg1$, the agreement with the logistic curve becomes particularly
good as the ambient temperature decreases. Note that this good agreement
corresponds to short relaxation times of the environmental variable,
$\gamma^{-1}=0.02\ll1$, which is in line with our theoretical assumptions.
The variation of the time scales $\gamma^{-1}$ and $\lambda^{-1}$
requires further research. Figure \ref{stlog2} confirms the validity
of the ergodic hypothesis. It can be seen that the numerically obtained
population density distribution fits the theoretical curve.

\section{Discussion}

On a fairly simple but conceptually important example, we have presented
the way of probabilistic reasoning and the mathematical structure
underlying the statistical approach in modelling evolutionary processes
in ecology and population dynamics. This example is rather a particular
implementation of the idea of the proposed statistical approach. In
this context, a natural question arises: What perspectives, in a wide
sense, does the statistical approach potentially provide in modelling
ecological processes, if any? Indeed, if mathematical ideas are formulated
in general terms, then it is necessary to clarify the further prospect
of their possible application.

Population growth models differ mainly in their phase space geometry,
which determines the feasible stationary solutions and the behaviour
of the system near these points. In statistical interpretation, this
implies a difference in transient processes to the most probable state
of the population. This distinction is central to the understanding
of ecological evolution in general. Having a range of mathematical
descriptions of transient processes provides an opportunity to implement
phase space geometry models and corresponding evolution that differ
in quality. In this context, let us consider the following concrete
example in order to answer, at least in part, the question posed.

Without over-complicating the evolutionary dynamics \prettyref{eq:StatLog},
we modify these equations on the assumption of non-Gaussian statistics
of the environment variable, while keeping the logistic $h$--function
\eqref{eq:log-h} unchanged. More precisely, we consider the equilibrium
(invariant) density of the form, 
\begin{equation}
\sigma_{\vartheta}^{+}(x,y)\propto\exp\left\{ -\vartheta^{-1}h(x)\right\} \times\exp\left\{ -\vartheta^{-1}h^{*}(y)\right\} ,\label{eq:nonGauss-y}
\end{equation}
where $h^{*}(y)$ is the environment $h$--function. Let the $\vartheta$-expression
be chosen in the form 
\[
\Theta^{*}(y,\vartheta)=\lambda\vartheta h^{*}\phantom{}'(y)-\gamma\left[\left(h^{*}\phantom{}'(y)\right)^{2}-\vartheta h^{*}\phantom{}''(y)\right].
\]
The proof that $\Theta^{*}(y,\vartheta)$ is indeed a $\vartheta$-expression
is done by direct calculation \citep{samoletov2021advanced}. Under
these assumptions, we arrive at the (stochastic) equations of motion,
\begin{eqnarray}
\dot{x} & = & \lambda xh^{*}\phantom{}'(y),\nonumber \\
\dot{y} & = & -\lambda\left[xh'(x)-\vartheta\right]-\gamma h^{*}\phantom{}'(y)+\sqrt{2\gamma\vartheta}\xi(t).\label{eq:nonGauss-Dyn}
\end{eqnarray}
The density $\sigma_{\vartheta}^{+}(x,y)$ \prettyref{eq:nonGauss-y}
is invariant for dynamics \prettyref{eq:nonGauss-Dyn}. The proof
is by direct calculation. In the case $h^{*}(y)=\nicefrac{y^{2}}{2}$,
$\vartheta$-expression and dynamic equations coincide with those
considered earlier. Note that if we put $\gamma\equiv0$, we obtain
a system of ordinary differential equations with the required invariant
density (easily checked by direct calculation). However, in this case
the ergodicity condition is problematic \citep{legoll2009non}.

The choice of admissible function $h^{*}\left(y\right)$ introduces
noticeable freedom in the modelling of transients. As an example,
let us consider two variants of the bimodal probability density function
$\sigma_{\vartheta}^{*}(y)$, symmetric and asymmetric. First, let
the symmetric density be defined by the expression 
\begin{equation}
h^{*}\phantom{}'(y)=y\left(y+\sqrt{m}\right)\left(y-\sqrt{m}\right),\label{eq:symm}
\end{equation}
where $m>0$ is a parameter. For the shortened ordinary differential
equation, $\dot{y}=-\gamma h^{*}\phantom{}'(y)$, the equilibria $y=\pm\sqrt{m}$
are stable, while $y=0$ is an unstable equilibrium.

To deepen the discussion, consider a hypothetical situation where
changes in the environment are associated not only with population
size, but also with changes in the phenotype of the population itself,
when an additional environmental resource becomes available to the
population and affects its growth and terminal size. The asymptotic
result would be a change in the value of the most probable population
density. We now focus on describing such a process. 
\begin{figure}[h]
\centering{}\includegraphics[scale=0.4]{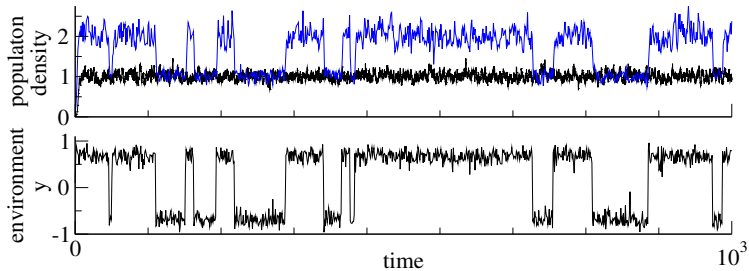}\caption{(colour online) Population density $x$ (top) and environment variable
$y$ (bottom) versus time for the symmetric bimodal equilibrium probability
density of variable y (equation\prettyref{eq:symm}) in two modes:
in the absence (black) and in the presence of functional population--environment
relationship (blue). The parameters used in the simulation are as
follows: $m=0.5,\vartheta=0.01,r=1,\gamma=50,\lambda=\sqrt{50}$.}
\end{figure}

To solve this problem, we assume that the functional relationship
between the population and the environment has the form $K\left[y\right]=1+H(y)$,
where $H(y)$ is a Heaviside function. Although this piecewise linear
relationship is speculative and oversimplified, it reflects an important
qualitative feature of the phenotype--environment relationship responsible
for the transition between two equilibrium population densities. It
would be more correct to relate $K$ to its own dynamics, but for
our illustrative purposes such a complication is unnecessary.

Qualitatively, this dependence can be explained as follows: the transition
of the environment from one steady state to another is associated
with a change in the phenotype of the population such that an additional
environmental resource becomes available to the population.

To test the intended perspective, we perform numerical simulations.
The parameters used in the simulation were chosen to remain close
to the logistic dynamics. Figure 3 shows the population density $x$
(top) simultaneously with the corresponding environmental variable
$y$ (bottom) versus time, in two different modes: in the absence
(black curve) and in the presence (blue curve) of a functional relationship
between population and environment $K[y]$. The random fluctuations
of the population density synchronised with the environment can be
clearly identified.

This behaviour provides the basis for the following important modification
of the bimodal density $\sigma_{\vartheta}^{*}(y)$. Namely, a conceptually
important example is the asymmetric bimodal probability density $\sigma_{\vartheta}^{*}(y)$,
defined by the expression 
\begin{equation}
h^{*}\phantom{}'(y)=Dy\left(y-a\right)\left(y-1\right).\label{eq:asymm}
\end{equation}
The detuning parameter $a$ in this expression satisfies $0<a<1$,
and the parameter $D>0$. Thus, for the shortened ordinary differential
equation, $\dot{y}=-\gamma h^{*}\phantom{}'(y)$, the equilibria $y=0,1$
are stable, while $y=a$ is an unstable equilibrium. Let there be
a functional relationship between the population and the environment
of the form $K\left[y\right]=1+H(y-a)$. This piecewise linear dependence
should, under certain conditions, determine the long-term transition
from the initial equilibrium population density to the new one. To
check this, we carry out numerical simulations. The parameters used
in the simulation are chosen to remain close to the logistic equation.
Figure 4 shows population density $x$ as a function of time in two
different modes: in the absence (black curve) and in the presence
(blue curve) of functional relationship $K[y]$ between population
and environment. For the selected asymmetric bimodal density parameters,
it is observed a clear transition to a new equilibrium state in the
population size. Mathematically, a backward transition to the initial
population size is possible, but the average waiting time is expected
to be long and may exceed the lifetime of the population \citep{stratonovich1963topics}.
We omit the discussion of the relevant mathematical details. 
\begin{figure}[h]
\centering{}\includegraphics[clip,scale=0.4]{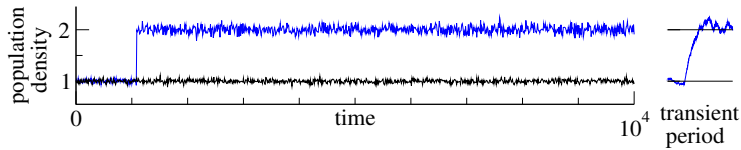}\caption{(colour online) Population density versus time for the equilibrium
asymmetric bimodal probability density of variable $y$ (equation
\eqref{eq:asymm}), in two different modes: in the absence (black)
and in the presence (blue) of functional population--environment
relationship. Note that the transition to a new population density
persists for a long time ("forever"). The
parameters used in the simulation are as follows: $D=4,a=0.25,\vartheta=0.001,r=1,\gamma=50,\lambda=\sqrt{50}$.}
\end{figure}

When we use the term ecological temperature, we refer to the parameter
$\vartheta$, which in turn is related to the key definition of the
$\vartheta$--expression in \prettyref{sec:Mathematical-formalism}.
For the presented statistical interpretation of the logistic equation,
the specific value of $\vartheta$ is immaterial. To give $\vartheta$
an appropriate value, one should choose a special device, an ecological
thermometer. The analogy with thermodynamics can be useful. However,
a discussion of this issue is beyond the scope of this article. In
an extended context, other parameters in addition to $\vartheta$
may be necessary, for example, when considering the processes of emigration
and immigration of the population.

Finally, we should make a note about the random term in the system
\prettyref{eq:nonGauss-Dyn}. If we put $\gamma\equiv0$, we obtain
a system of ordinary differential equations for which the density
\prettyref{eq:nonGauss-y} is invariant (this is easy to check). However,
these equations of motion are not ergodic for the given density (note
that there exists the integral of motion, $I=h(x)-\vartheta\ln x+h^{*}(y)=const$),
so the dynamical scenarios discussed earlier have no place.

\section{Closing remarks}

This paper introduces a mathematical scheme for the statistical interpretation
of growth models, by using the logistic equation as an example. The
statistical interpretation of the Gompertz model \citep{Winsor1932,tjorve2017use},
which is the second most common growth model, can be treated similarly.
Under the same assumptions as for the logistic model, first the Gompertz
$h$--function is determined, that is,
\[
h(x)=\frac{1}{2K}\left(\ln\frac{x}{K}\right)^{2}
\]
Thus, the stochastic equations of the Gompertz growth model take the
form:
\begin{eqnarray*}
\dot{x} & = & \lambda xy,\\
\dot{y} & = & -\lambda\left[\ln\left(\frac{x}{K}\right)-\vartheta\right]-\gamma y+\sqrt{2\gamma\vartheta}\xi(t).
\end{eqnarray*}
The Gompertz model is often used to describe the growth of plants
and animals, bacteria and cancer cells, and other processes (\emph{e.g.},
\citep{Chaplain2014mathematical,vasiev2016aging}).

Note that the $h$--functions of both the Gompertz and logistic growth
models have only one parameter, the carrying capacity, $K$. However,
when fitting experimental data, a model that offers greater flexibility
in growth dynamics may be preferable. This is not an issue in the
context of this paper. 

Anyway, suppose an experimental observation yields
a normalised histogram of the population size. For example, this can be
extracted from a time series of population size changes that occur
under ecological equilibrium conditions. From this, we can calculate
the $h$--function corresponding to the histogram, which can be approximated
by a curve with two or more parameters. Using this $h$--function,
we can write down the equations for population growth. Inevitably,
experimental data sets are limited. Thus, it is essential to evaluate the 
predictive ability of the derived dynamic equations 
to drive additional model fitting. 

\section*{Appendix \label{sec:Appendix-A} }

We are interested in the dynamic equations,
\begin{equation}
\dot{x}=g(x,y),\quad\dot{y}=g^{*}(x,y),\label{eq:a1}
\end{equation}
provided the functional equation of the dynamic principle,
\begin{equation}
h'\left(x\right)g(x,y)+yg^{*}(x,y)=c_{1}y+c_{2}(y^{2}-\vartheta),\label{eq:a2}
\end{equation}
is satisfied, and the probability density $\sigma_{\vartheta}^{+}(x,y)$\eqref{eq:probdensity}
is invariant for the dynamics \eqref{eq:a1}. The problem is to find
suitable functions $g(x,y),g^{*}(x,y)$ such that the equation \eqref{eq:a2}
is valid. 

Let us divide the procedure for finding the functions $g(x,y),g^{*}(x,y)$
into two steps. First we represent the function $g^{*}(x,y)$ as the
sum of two functions,
\begin{equation}
g^{*}(x,y)=g_{0}^{*}(x)+g_{1}^{*}(y).\label{eq:split-g}
\end{equation}
Then we have to satisfy two equations separately:
\begin{eqnarray}
h'\left(x\right)g(x,y)+yg_{0}^{*}(x) & = & c_{1}y,\label{eq:g0}\\
yg_{1}^{*}(y) & = & c_{2}(y^{2}-\vartheta).\label{eq:g1}
\end{eqnarray}
The last equation does not have a deterministic solution. So we have
to consider random functions and replace equation \eqref{eq:g1} with
\begin{equation}
\mathbb{E}_{\omega}\left\{ yg_{1}^{*}(y)-c_{2}(y^{2}-\vartheta)\right\} =0.\label{eq:g1rand}
\end{equation}

We start by solving equation \eqref{eq:g0}. Substituting 
\[
g(x,y)=y\varphi(x),
\]
where $\varphi(x)$ is a function, in equation \eqref{eq:g0}, we
obtain the following equation,
\[
\varphi(x)h'\left(x\right)-c_{1}=-g_{0}^{*}(x).
\]
Note that the expression $c_{1}y$ does not define a value for the
ecological temperature. Therefore, we must additionally require that
the probability density $\sigma_{\vartheta}^{+}(x,y)$\eqref{eq:probdensity}
is invariant for the deterministic dynamics,
\begin{equation}
\dot{x}=y\varphi(x),\quad\dot{y}=-\left[\varphi(x)h'\left(x\right)-c_{1}\right],\label{eq:detsys0}
\end{equation}
which is equivalent to the requirement: $\sigma_{\vartheta}^{+}(x,y)$
is a stationary solution of the Liouville equation corresponding to
the dynamical system \eqref{eq:detsys0}, that is,
\begin{equation}
-\frac{\partial}{\partial x}\left[y\varphi(x)\sigma_{\vartheta}^{+}(x,y)\right]-\frac{\partial}{\partial y}\left[-\left(\varphi(x)h'\left(x\right)-c_{1}\right)\sigma_{\vartheta}^{+}(x,y)\right]=0.\label{eq:Liouville}
\end{equation}
After computing the derivatives, it can be deduced that equation \eqref{eq:Liouville}
is an identity only if
\[
c_{1}=\varphi'(x)\vartheta.
\]
This leaves uncertainty in the function $\varphi(x)$. Assuming $c_{1}$
to be a constant, we derive the following form of the function $\varphi(x)$:
\[
\varphi(x)=\lambda x,
\]
where $\lambda$ is a constant coefficient. Therefore, we have arrived
at the expression
\begin{equation}
g_{0}^{*}(x)=-\lambda\left[xh'\left(x\right)-\vartheta\right].\label{eq:g0-x}
\end{equation}

Next, we need to determine the function $g_{1}^{*}(y)$. By substituting
\[
g_{1}^{*}(y)=c_{2}y+\psi
\]
into the equation \eqref{eq:g1rand}, we arrive at an equation that
has no deterministic solution. Treating $\psi$ as a random variable
requires us to examine the following equation,
\begin{equation}
\mathbb{E}_{\omega}\left\{ y\psi\right\} =-c_{2}\vartheta.\label{eq:E-theta-y}
\end{equation}

The equation,
\begin{equation}
\dot{y}=g^{*}(x,y)=g_{0}^{*}(x)+c_{2}y+\psi,\label{eq:y_psi}
\end{equation}
shows that the variable $y$ is functionally dependent on $\psi$.
Thus, by meeting certain assumptions about the random variable $\psi$,
Novikov's formula \citep{klyatskin2005dynamics} can be applied to
calculate the mathematical expectation $\mathbb{E}_{\omega}\left\{ y\psi\right\} $.
The most straightforward way to implement $\psi$ is to use a white
noise process as $\psi$. Indeed, let us set
\[
\psi=a\xi(t),
\]
where $\xi(t)$ is the standard Gaussian white noise, $\mathbb{E}_{\omega}\left\{ \xi(t)\right\} =0$,
$\:\mathbb{E}_{\omega}\left\{ \xi(t)\xi(t')\right\} =\delta(t-t')$,
and $a$ is a constant coefficient. The equation \eqref{eq:y_psi}
takes the form of a stochastic differential equation with additive
noise. If that is the case, Novikov's formula \citep{klyatskin2005dynamics}
takes a particularly simple form:
\[
\mathbb{E}_{\omega}\left\{ y\psi\right\} =a\mathbb{E}_{\omega}\left\{ y\xi(t)\right\} =\frac{1}{2}a^{2},
\]
where the expression of the functional derivative below \citep{klyatskin2005dynamics}
has been taken into account,
\[
\frac{\delta y(t)}{\delta\xi(t)}=a.
\]
To calculate this functional derivative, we should first write equation \eqref{eq:y_psi}) in integral form and then take into account the fact that $y(t')$ is functionally independent of $\xi(t)$ at $t'<t$.

Thus, we arrive at the relationship,
\[
\frac{1}{2}a^{2}=-c_{2}\vartheta,
\]
that means $c_{2}<0$ . 

A particular solution of the dynamic principle equation as the stochastic
equations of motion \eqref{eq:StochPopDyn} is obtained by denoting
$c_{2}=-\gamma$ where $\gamma>0$, and collecting and substituting
all results into equations \eqref{eq:a1}, that is,
\begin{eqnarray}
	\dot{x} & = & \lambda xy,\nonumber \\
	\dot{y} & = & -\lambda\left[xh'(x)-\vartheta\right]-\gamma y +\sqrt{2\gamma\vartheta}\xi(t)\nonumber.
\end{eqnarray}

\section*{Acknowledgement}

This work has been supported by the EPSRC grant EP/S033211/1.

\bibliographystyle{elsarticle-num}

\begin{thebibliography}{30}
	
	\bibitem{Verhulst1838}
	P.-F. Verhulst, Notice sur la loi que la population suit dans son
	accroissement, Corresp. Math. Phys. 10 (1838) 113--126.
	
	\bibitem{pearl1924curve}
	R.~Pearl, The curve of population growth, Proceedings of the American
	Philosophical Society 63~(1) (1924) 10--17.
	
	\bibitem{kingsland1982}
	S.~Kingsland, The refractory model: The logistic curve and the history of
	population ecology, The Quarterly Review of Biology 57~(1) (1982) 29--52.
	
	\bibitem{kot2001elements}
	M.~Kot, Elements of mathematical ecology, Cambridge University Press, 2001.
	
	\bibitem{Murray2002}
	J.~D. Murray, Mathematical Biology I. An Introduction, Vol.~17, Springer, New
	York, 2002.
	
	\bibitem{volterra1939calculus}
	V.~Volterra, Calculus of variations and the logistic curve, Human Biology
	11~(2) (1939) 173--178.
	
	\bibitem{leitmann1972minimum}
	G.~Leitmann, A minimum principle for a population equation, Journal of
	Optimization Theory and Applications 9~(2) (1972) 155--156.
	
	\bibitem{gatto1988functional}
	M.~Gatto, S.~Muratori, S.~Rinaldi, A functional interpretation of the logistic
	equation, Ecological Modelling 42~(2) (1988) 155--159.
	
	\bibitem{webb1995hamilton}
	J.~N. Webb, Hamilton's variational principle and ecological models, Ecological
	modelling 80~(1) (1995) 35--40.
	
	\bibitem{pawlowski2006dynamic}
	C.~W. Pawlowski, Dynamic landscapes, stability and ecological modeling, Acta
	Biotheoretica 54~(1) (2006) 43--53.
	
	\bibitem{wilhelm2000goal}
	T.~Wilhelm, R.~Br{\"u}ggemann, Goal functions for the development of natural
	systems, Ecological Modelling 132~(3) (2000) 231--246.
	
	\bibitem{may2019stability}
	R.~M. May, Stability and complexity in model ecosystems, Princeton university
	press, 2019.
	
	\bibitem{SamoletovVasiev2017}
	A.~Samoletov, B.~Vasiev, Dynamic principle for ensemble control tools, J. Chem.
	Phys. 147~(20) (2017) 204106.
	
	\bibitem{samoletov2021advanced}
	A.~Samoletov, B.~Vasiev, Advanced selection of ensemble control tools, Journal
	of Physics: Conference Series 2090 (2021) 012059.
	
	\bibitem{SamoletovDettmannChaplain2007}
	A.~Samoletov, C.~Dettmann, M.~Chaplain, Thermostats for "slow" configurational
	modes, J. Stat. Phys. 128~(6) (2007) 1321--1336.
	
	\bibitem{SamoletovDettmannChaplain2010}
	A.~Samoletov, C.~Dettmann, M.~Chaplain, Notes on configurational thermostat
	schemes, J. Chem. Phys. 132~(24) (2010) 246101.
	
	\bibitem{hastings2004transients}
	A.~Hastings, Transients: the key to long-term ecological understanding?, Trends
	in ecology \& evolution 19~(1) (2004) 39--45.
	
	\bibitem{Hastings2018}
	A.~Hastings, K.~C. Abbott, K.~Cuddington, T.~Francis, G.~Gellner, Y.-C. Lai,
	A.~Morozov, S.~Petrovskii, K.~Scranton, M.~L. Zeeman, Transient phenomena in
	ecology, Science 361~(6406) (2018) eaat6412.
	
	\bibitem{legoll2009non}
	F.~Legoll, M.~Luskin, R.~Moeckel, Non-ergodicity of {N}os{\'e}--{H}oover
	dynamics, Nonlinearity 22~(7) (2009) 1673.
	
	\bibitem{khinchin1949mathematical}
	A.~Y. Khinchin, Mathematical foundations of statistical mechanics, Dover
	Publications, Inc., New York, 1949.
	
	\bibitem{gardiner2009stochastic}
	C.~W. Gardiner, Stochastic Methods: A Handbook for the Natural and Social
	Sciences, Springer, Berlin, 2009.
	
	\bibitem{klyatskin2005dynamics}
	V.~I. Klyatskin, Dynamics of stochastic systems, Elsevier, 2005.
	
	\bibitem{kramers1940}
	H.~A. Kramers, Brownian motion in a field of force and the diffusion model of
	chemical reactions, Physica 7~(4) (1940) 284--304.
	
	\bibitem{Samoletov1999}
	A.~A. Samoletov, A remark on the {K}ramers problem, J. Stat. Phys. 96~(5-6)
	(1999) 1351--1357.
	
	\bibitem{mel1991kramers}
	V.~I. Mel'nikov, The {K}ramers problem: Fifty years of development, Physics
	Reports 209~(1-2) (1991) 1--71.
	
	\bibitem{stratonovich1963topics}
	R.~L. Stratonovich, Topics in the theory of random noise, Vol.~1, Gordon and
	Breach, 1963.
	
	\bibitem{Winsor1932}
	C.~P. Winsor, The Gompertz curve as a growth curve, Proc. Natl. Acad. Sci.
	U.S.A. 18~(1) (1932) 1.
	
	\bibitem{tjorve2017use}
	K.~M. Tj{\o}rve, E.~Tj{\o}rve, The use of Gompertz models in growth analyses,
	and new Gompertz-model approach: An addition to the unified-richards family,
	PloS one 12~(6) (2017) e0178691.
	
	\bibitem{Chaplain2014mathematical}
	H.~Enderling, M.~Chaplain, Mathematical modeling of tumor growth and
	treatment, Current Pharmaceutical Design 20~(30) (2014) 4934--4940.
	
	\bibitem{vasiev2016aging}
	D.~Avraam, S.~Arnold, O.~Vasieva, B.~Vasiev, On the heterogeneity of human
	populations as reflected by mortality dynamics, Aging (Albany NY) 8~(11)
	(2016) 3045.
	
\end{thebibliography}

\end{document}